\documentclass[12pt]{amsart}
\usepackage{amssymb,amsmath,amsthm,mathtools,amsfonts,times,hyperref,mathrsfs,multirow,tabularx}
\usepackage{pgfplots}
\pgfplotsset{compat=1.15}
\usepackage{mathrsfs,mathdots}
\usetikzlibrary{arrows}
\usepackage{enumerate}
\usepackage{graphicx}
\usepackage{array}
\usepackage{ytableau}
\usepackage{pstricks,pst-node,graphicx}
\usepackage{tikz,varwidth}
\usepackage{setspace}
\usepackage{float}
\usetikzlibrary{calc}
\usepackage[super]{nth}
\usepackage{tikz}
\usepackage{extarrows}
\usepackage{pgfplots}
\pgfplotsset{compat=1.15}
\usetikzlibrary{arrows}

\newtheorem{theorem}{Theorem}[section]
\newtheorem{lemma}[theorem]{Lemma}
\newtheorem{cor}[theorem]{Corollary}
\newtheorem{conj}[theorem]{Conjecture}
\newtheorem{prop}[theorem]{Proposition}
\theoremstyle{definition}

\newtheorem{example}[theorem]{Example}

\theoremstyle{remark}
\newtheorem{remark}{Remark}
\numberwithin{equation}{section}

\allowdisplaybreaks

\newcommand{\at}[2][]{#1|_{#2}}

\newcommand\numberthis{\addtocounter{equation}{1}\tag{\theequation}}
\newcommand\restr[2]{{
  \left.\kern-\nulldelimiterspace 
  #1 
  \littletaller 
  \right|_{#2} 
  }}
\newcommand{\littletaller}{\mathchoice{\vphantom{\big|}}{}{}{}}  

\begin{document}

\title[Finite analogs of partition inequalities related to hook length two]{Finite analogs of partition bias related to hook length two and a variant of Sylvester's map}

\author{Alexander Berkovich}
\address{Department of Mathematics, University of Florida, Gainesville
FL 32611, USA}
\email{alexb@ufl.edu}
\author{Aritram Dhar}
\address{Department of Mathematics, University of Florida, Gainesville
FL 32611, USA}
\email{aritramdhar@ufl.edu}

\dedicatory{Dedicated to Krishnaswami Alladi with gratitude and admiration}

\date{\today}



\begin{abstract}
In this paper, we count the total number of hooks of length two in all odd partitions of $n$ and all distinct partitions of $n$ with a bound on the largest part of the partitions. We generalize inequalities of Ballantine, Burson, Craig, Folsom and Wen by showing there is a bias in the number of hooks of length two in all odd partitions over all distinct partitions of $n$ in presence of a bound on the largest part. To establish such a bias, we use a variant of Sylvester's map. Then, we conjecture a similar finite bias for a weighted count of hooks of length two and prove it when we remove the bound on the largest part.   
\end{abstract}

\maketitle

\section{Introduction}\label{s1}
A \textit{partition} $\pi$ is a non-increasing finite sequence $\pi = (\lambda_1,\lambda_2,\lambda_3,\dots)$ of positive integers. The elements $\lambda_i$ appearing in the sequence $\pi$ are called the \textit{parts} of $\pi$. The number of parts of $\pi$ is denoted by $\nu(\pi)$ and the largest part of $\pi$ is denoted by $l(\pi)$. The sum of all the parts of $\pi$ is called the \textit{size} of $\pi$ and is denoted by $|\pi|$. The \textit{alternating sum of parts} of $\pi$, denoted by $\gamma(\pi)$, is defined as
\begin{align*}
\gamma(\pi) = \lambda_1-\lambda_2+\lambda_3-\ldots.    
\end{align*}
We say $\pi$ is a partition of $n$ if its size is equal to $n$. The empty sequence $\emptyset$ is considered as the unique partition of zero. Let $\mathcal{P}$ denote the set of all partitions.\\\par
\begin{figure}[H]
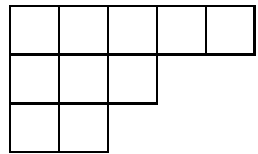

\ytableausetup{centertableaux}
\ytableaushort
{{}{}{}{}{},{}{}{},{}{}}
* {5,3,2}
\par\quad\\\par\quad\\
\caption{The Young diagram of the partition $\pi = (5,3,2)$.}
\label{fig1}
\end{figure}

 The \textit{Young diagram} of a partition $\pi$ is a way of representing $\pi$ graphically where the parts $\lambda_i$ of $\pi$, called \textit{cells}, are arranged in left-justified rows with $\lambda_i$ cells in the $i$-th row.  See Figure \ref{fig1} below for an example of a Young diagram of a partition.\\
\par Let us now recall \textit{Sylvester's bijection} \cite{B99,P06}. Define $\mathcal{O}_{n}$ (resp. $\mathcal{D}_{n}$) to be the set of all partitions $\pi_o$ (resp. $\pi_d$) of $n$ into odd (resp. distinct) parts.\\\par Consider the map $\psi : \mathcal{O}_{n}\rightarrow\mathcal{D}_{n}$ and let $\pi_{o}\in\mathcal{O}_{n}$. We arrange all odd parts of $\pi_{o}$ symmetrically row-by-row and then read the number of cells along the \textit{center hook} starting from the last row following the \textit{spine} to the top and then to the right. Then, read the number of cells along the hook to the left of the spine from bottom to top and left followed by the number of cells along the hook to the right of the spine from bottom to top and right and so on. Finally, we get a partition $\pi_{d}$ of $n$ into distinct parts, i.e., $\pi_{d}\in\mathcal{D}_n$.\\\par Clearly, $\psi$ is invertible. Consider $\psi^{-1} := \sigma : \mathcal{D}_n\rightarrow\mathcal{O}_n$ and let $\pi_{d} = (\lambda_1,\lambda_2,\lambda_3,\ldots)\in\mathcal{D}_n$. Now, consider the alternating sum of parts $\gamma(\pi_d)$ of $\pi_d$. We then construct a center hook consisting of $\gamma(\pi_d)$ cells placed vertically from bottom to top to form a spine and $\lambda_1-\gamma(\pi_d)$ cells placed horizontally to the right of the topmost cell of the spine. Then, symmetrically insert $\lambda_1-\gamma(\pi_d)$ cells horizontally to the left of the topmost cell of the spine. Then, add $\lambda_2-\lambda_1+\gamma(\pi_d)$ cells vertically from top to bottom below the rightmost cell of the horizontal piece to the left of the spine. Then, symmetrically insert $\lambda_2-\lambda_1+\gamma(\pi_d)$ cells vertically from top to bottom starting from the cell in second row and in second column to the right of the spine. Now, add $\lambda_3-\lambda_2+\lambda_1-\gamma(\pi_d)$ cells to the right of the cell in second row and in second column to the right of the spine. Proceeding in this fashion, after exhausting all the parts of $\pi_d$, we read the number of cells in the newly constructed symmetric diagram row-by-row from top to bottom to get parts of a partition $\pi_o$ of $n$ into odd parts, i.e., $\pi_{o}\in\mathcal{O}_n$.\\\par
\begin{figure}[H]
\centering
\definecolor{qqwuqq}{rgb}{0,0.39215686274509803,0}
\definecolor{ffqqqq}{rgb}{1,0,0}
\begin{tikzpicture}[line cap=round,line join=round,>=triangle 45,x=0.8cm,y=0.8cm]
\clip(-5.1,-5.6) rectangle (20.0,-1.0);
\draw [line width=1.5pt] (0.8310151627509912,-1.2771841301361357)-- (7.83101516275099,-1.2771841301361357);
\draw [line width=1.5pt] (0.8310151627509912,-2.2771841301361357)-- (7.83101516275099,-2.2771841301361357);
\draw [line width=1.5pt] (1.8310151627509912,-3.2771841301361357)-- (6.83101516275099,-3.2771841301361357);
\draw [line width=1.5pt] (2.831015162750991,-4.277184130136137)-- (5.83101516275099,-4.277184130136137);
\draw [line width=1.5pt] (2.831015162750991,-5.277184130136137)-- (5.83101516275099,-5.277184130136137);
\draw [line width=1.5pt] (0.8310151627509912,-1.2771841301361357)-- (0.8310151627509912,-2.2771841301361357);
\draw [line width=1.5pt] (1.8310151627509912,-1.2771841301361357)-- (1.8310151627509912,-3.2771841301361357);
\draw [line width=1.5pt] (2.831015162750991,-1.2771841301361357)-- (2.831015162750991,-5.277184130136137);
\draw [line width=1.5pt] (3.831015162750991,-1.2771841301361357)-- (3.831015162750991,-5.277184130136137);
\draw [line width=1.5pt] (4.83101516275099,-1.2771841301361357)-- (4.83101516275099,-5.277184130136137);
\draw [line width=1.5pt] (5.83101516275099,-1.2771841301361357)-- (5.83101516275099,-5.277184130136137);
\draw [line width=1.5pt] (6.83101516275099,-1.2771841301361357)-- (6.83101516275099,-3.2771841301361357);
\draw [line width=1.5pt] (7.83101516275099,-1.2771841301361357)-- (7.83101516275099,-2.2771841301361357);
\draw [line width=1.5pt,color=ffqqqq] (4.35101516275099,-1.7971841301361358)-- (4.35101516275099,-4.777184130136137);
\draw [line width=1.5pt,color=ffqqqq] (4.35101516275099,-1.7971841301361358)-- (7.35101516275099,-1.7771841301361357);
\draw [line width=1.5pt,color=ffqqqq] (5.35101516275099,-2.7771841301361357)-- (5.33101516275099,-4.777184130136137);
\draw [line width=1.5pt,color=ffqqqq] (5.35101516275099,-2.7771841301361357)-- (6.35101516275099,-2.7771841301361357);
\draw [line width=1.5pt,color=qqwuqq] (1.3710151627509912,-1.7771841301361357)-- (3.331015162750991,-1.7971841301361358);
\draw [line width=1.5pt,color=qqwuqq] (3.331015162750991,-1.7971841301361358)-- (3.331015162750991,-4.777184130136137);
\begin{scriptsize}
\draw [fill=ffqqqq] (4.35101516275099,-4.777184130136137) circle (2.5pt);
\draw [fill=ffqqqq] (4.35101516275099,-1.7971841301361358) circle (2.5pt);
\draw [fill=ffqqqq] (7.35101516275099,-1.7771841301361357) circle (2.5pt);
\draw [fill=ffqqqq] (5.33101516275099,-4.777184130136137) circle (2.5pt);
\draw [fill=ffqqqq] (5.35101516275099,-2.7771841301361357) circle (2.5pt);
\draw [fill=ffqqqq] (6.35101516275099,-2.7771841301361357) circle (2.5pt);
\draw [fill=qqwuqq] (3.331015162750991,-4.777184130136137) circle (2.5pt);
\draw [fill=qqwuqq] (3.331015162750991,-1.7971841301361358) circle (2.5pt);
\draw [fill=qqwuqq] (1.3710151627509912,-1.7771841301361357) circle (2.5pt);
\draw [fill=qqwuqq] (2.331015162750991,-2.7571841301361357) circle (2.5pt);
\end{scriptsize}
\end{tikzpicture}
\caption{An example of Sylvester's bijection: $\pi_o = (7,5,3,3)$ (read row by row). $\pi_d = \psi(\pi_o) = (7,6,4,1)$. The red lines are the odd-indexed parts of $\pi_d$ and the green lines are the even-indexed parts of $\pi_d$. It is clear that $\gamma(\pi_d) = 4$.}
\label{fig2}
\end{figure}
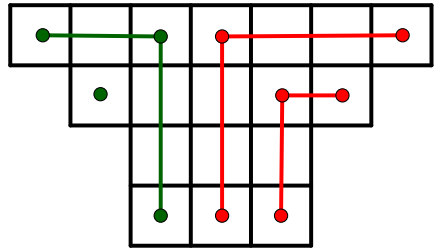

 See Figure \ref{fig2} below for an example illustrating Sylvester's map $\psi$.\\
\par The conjugate of a partition $\pi$ is the partition $\pi^{\prime} = (\lambda^{\prime}_1,\lambda^{\prime}_2,\lambda^{\prime}_3,\dots)$ whose Young diagram has the columns of $\pi$ as rows. Each cell in a Young diagram of $\pi$ may be labeled with a number called \textit{hook length}, which is one plus the number of cells directly to the right and directly below the chosen cell. For a cell in the $i$-th row and $j$-th column of the Young diagram of a partition $\pi$, its \textit{hook length} is defined as $h(i,j) = \lambda_i + \lambda^{\prime}_{j} - i - j + 1$. Let $\mathcal{H}(\pi)$ denote the multi-set of all hook lengths of $\pi$. See Figure \ref{fig3} below for an example of a Young diagram of a partition along with its hook lengths.\\
\begin{figure}[H]
\ytableausetup{centertableaux}
\ytableaushort
{{12}{10}{8}{7}{6}{4}{2}{1},{9}{7}{5}{4}{3}{1},{7}{5}{3}{2}{1},{3}{1},{1}}
* {8,6,5,2,1}
\par\quad\\\par\quad\\
\caption{The Young diagram of the partition $\pi = (8,6,5,2,1)$ with its hook lengths.}
\label{fig3}
\end{figure}
\par In $2006$, Nekrasov and Okounkov \cite{NO06} discovered a formula for arbitrary powers of Euler’s infinite product in terms of hook lengths which can be stated as follows.\\
\begin{theorem}\label{thm11}
For any complex number $z$, we have\\
\begin{align*}
\sum\limits_{\pi\in\mathcal{P}}x^{|\pi|}\prod_{h\in\mathcal{H}(\pi)}\left(1-\dfrac{z}{h^2}\right) = \prod_{k=1}^{\infty}(1-x^k)^{z-1}.\numberthis\label{eq11}\\    
\end{align*}    
\end{theorem}
\par Fix any positive integer $t$. Let $a_t(n)$ (resp. $b_t(n)$) denote the total number of hooks of length $t$ in all partitions $\pi_o$ (resp. $\pi_d$) of $n$ into odd parts (resp. distinct parts).
\\\par For any partition $\pi$, a cell in its Young diagram has hook length $1$ if and only if it is at the end of a row and there is no cell directly below it. Thus, in a partition $\pi$, the number of hooks of length $1$ is equal to the number of different part sizes in $\pi$. The following result was conjectured by Beck \cite{O17} and proved analytically by Andrews \cite{A17}.\\
\begin{theorem} (Andrews \cite{A17})\label{thm12}
The difference between the total number of parts in all distinct partitions of $n$ and the total number of different part sizes in all odd partitions of $n$ equals $c(n)$, the number of partitions of $n$ with exactly one part occurring three times while all other parts occur only once.\\    
\end{theorem}
\begin{cor} (Ballantine et al. \cite[Corollary $1.6$]{BBCFW23})\label{cor13}
For all $n\ge 0$,\\
\begin{align*}
b_1(n) - a_1(n) = c(n).\numberthis\label{eq12}\\    
\end{align*}    
\end{cor}
\par From \eqref{eq12}, it is clear that $b_1(n)\ge a_1(n)$. Ballantine et al. \cite{BBCFW23} then stated the following conjecture which shows that the bias reverses for $t\ge 2$.\\
\begin{conj} (Ballantine et al. \cite[Conjecture $1.7$]{BBCFW23})\label{conj14}
\quad\\
\begin{enumerate}[(i)]
\item For every integer $t\ge 2$, there exists an integer $N_t$ such that for all $n > N_t$, we have $a_t(n)\ge b_t(n)$. Moreover, we conjecture the following values of $N_t$ for $2\le t\le 10$:\\
\begin{center}
\begin{tabularx}{.8\textwidth} { 
| >{\centering\arraybackslash}X 
| >{\centering\arraybackslash}X 
| >{\centering\arraybackslash}X
| >{\centering\arraybackslash}X
| >{\centering\arraybackslash}X
| >{\centering\arraybackslash}X
| >{\centering\arraybackslash}X
| >{\centering\arraybackslash}X
| >{\centering\arraybackslash}X 
| >{\centering\arraybackslash}X | }
\hline
$t$ & $2$ & $3$ & $4$ & $5$ & $6$ & $7$ & $8$ & $9$ & $10$\\
\hline
$N_t$ & $0$ & $7$ & $8$ & $18$ & $16$ & $34$ & $34$ & $56$ & $59$\\
\hline
\end{tabularx}\\\quad\\\quad\\
\end{center}   
\item For every integer $t\ge 2$, we have that $a_t(n)-b_t(n)\rightarrow\infty$ as $n\rightarrow\infty$.\\ 
\end{enumerate}
\end{conj}
Conjecture \ref{conj14} was recently proved by Craig, Dawsey and Han (see \cite[Theorem $1.2$]{CDH23}).\\\par Here and throughout the rest of the paper, for any two power series $A(q) := \sum\limits_{n\ge 0}a(n)q^n$ and $B(q) := \sum\limits_{n\ge 0}b(n)q^n$, $A(q)\succeq B(q)$ is equivalent to saying that $a(n)\ge b(n)$ for all $n\ge 0$.\\\par Let $L,m,n$ be non-negative integers. We now recall some notations from the theory of $q$-series that can be found in \cite{A98}. To this end, we define the conventional $q$-Pochhammer symbol as
\begin{align*}
(a)_L = (a;q)_L &:= \prod_{k=0}^{L-1}(1-aq^k),\\
(a)_{\infty} = (a;q)_{\infty} &:= \lim_{L\rightarrow \infty}(a)_L\,\,\text{where}\,\,\lvert q\rvert<1.
\end{align*}\\
We define the $q$-binomial (Gaussian) coefficient as\\
\begin{align*}
\left[\begin{matrix}m\\n\end{matrix}\right]_q := \Bigg\{\begin{array}{lr}
\dfrac{(q)_m}{(q)_n(q)_{m-n}}\quad\text{for } m\ge n\ge 0,\\
0\qquad\qquad\quad\text{otherwise}.\end{array}
\end{align*}\\
For $m, n\ge 0$, $\left[\begin{matrix}m+n\\n\end{matrix}\right]_q$ is the generating function for partitions into at most $n$ parts each of size at most $m$ (see \cite[Chapter $3$]{BBCFW23}).\\\par Ballantine et al. \cite{BBCFW23} derived the generating functions of $a_2(n)$ and  $b_2(n)$ which are as follows.\\
\begin{prop} (Ballantine et al. \cite[Proposition $3.1$]{BBCFW23})\label{prop15}
We have\\
\begin{align*}
\sum\limits_{n\ge 0}a_2(n)q^n = \dfrac{1}{(q;q^2)_{\infty}}\left(q^2+\sum\limits_{n\ge 2}\left(q^{2n-1}+q^{2(2n-1)}\right)\right).\numberthis\label{eq13}\\    
\end{align*}
\end{prop}
\begin{prop} (Ballantine et al. \cite[Proposition $3.2$]{BBCFW23})\label{prop16}
We have\\
\begin{align*}
\sum\limits_{n\ge 0}b_2(n)q^n = \dfrac{q^2}{1-q}(-q^2;q)_{\infty}.\numberthis\label{eq14}\\   
\end{align*}
\end{prop}
They showed that $a_2(n)\ge b_2(n)$ for all $n\ge 0$ using both $q$-series manipulations \cite[eq. ($3.5$)]{BBCFW23} and a combinatorial interpretation \cite[Proposition $3.3$]{BBCFW23}.\\\par Ballantine et al. \cite{BBCFW23} also gave the generating functions of $a_3(n)$ and $b_3(n)$ which are as follows.\\
\begin{prop} (Ballantine et al. \cite[Proposition $4.1$]{BBCFW23})\label{prop17}
We have\\
\begin{align*}
\sum\limits_{n\ge 0}a_3(n)q^n = (-q^3;q)_{\infty}\dfrac{q^3(1+q^3)}{1-q^2}+(-q;q)_{\infty}\left(\dfrac{q^6}{1-q^4}+\dfrac{q^3}{1-q^6}\right).\numberthis\label{eq15}\\ 
\end{align*}
\end{prop}
\begin{prop} (Ballantine et al. \cite[Proposition $4.2$]{BBCFW23})\label{prop18}
We have\\
\begin{align*}
\sum\limits_{n\ge 0}b_3(n)q^n = (-q;q)_{\infty}\sum\limits_{m\ge 2}\dfrac{q^m}{1+q^m}-\dfrac{q^2}{1-q^2}(-q^3;q)_{\infty}.\numberthis\label{eq16}\\
\end{align*}
\end{prop}
They also showed that $a_3(n)\ge b_3(n)$ for all $n\ge 7$ using $q$-series manipulations \cite[Proposition $4.3$]{BBCFW23}.\\\par For any positive integer $L$, we now define $a_t(L,n)$ (resp. $b_t(L,n)$) to be the total number of hooks of length $t$ in all partitions $\pi_o$ (resp. $\pi_d$) of $n$ into odd (resp. distinct) parts with $l(\pi_o)\le 2L-1$ (resp. $l(\pi_d)\le L$).\\\par Let $a_2(L,m,n)$ (resp. $b_2(L,m,n)$) denote the number of partitions $\pi_o$ (resp. $\pi_d$) of $n$ into odd (resp. distinct) parts with $l(\pi_o)\le 2L-1$ (resp. $l(\pi_d)\le L$) having $m$ hooks of length $2$.\\
\begin{example}
We compute $b_2(5,7) = 4$ and $b_3(5,7) = 4$. The partitions of $7$ into distinct parts each less than or equal to $5$ are $(5,2)$, $(4,3)$, and $(4,2,1)$. In Figure \ref{fig4}, we show their Young diagrams with labeled hook lengths.\\ 
\begin{figure}[H]
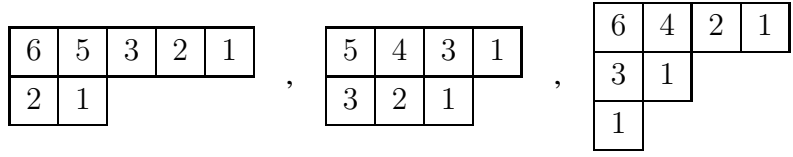

\ytableausetup{centertableaux}
\ytableaushort
{{6}{5}{3}{2}{1},{2}{1}}
* {5,2}
\quad,\quad
\ytableaushort
{{5}{4}{3}{1},{3}{2}{1}}
* {4,3}
\quad,\quad
\ytableaushort
{{6}{4}{2}{1},{3}{1},{1}}
* {4,2,1}
\par\quad\\\par\quad\\
\caption{The distinct partitions $\pi_d$ of $n=7$ having $l(\pi_d)\le 5$ and their hook lengths.}
\label{fig4}
\end{figure}    
\end{example}
\par In this paper, we consider finite analogs of Proposition \ref{prop15} and \ref{prop16} and consequently derive the generating functions of $a_2(L,n)$ and $b_2(L,n)$ in section \ref{s2}, namely, Theorems \ref{thm21} and \ref{thm22} respectively. We also show that the bias of the total number of hooks of length $2$ for odd partitions over distinct partitions remains the same in the finite case, i.e., $a_2(L,n)\ge b_2(L,n)$ for all $n\ge 0$, as stated in Theorem \ref{thm23}.\\\par The rest of the paper is organized as follows. In Section \ref{s2}, we present the statements of the main results and also state Conjecture \ref{conj26}. In Section \ref{s3}, we establish proofs of Theorems \ref{thm21} - \ref{thm25}, Corollary \ref{cor27} - \ref{cor28}, and Theorem \ref{thm29}. Finally, in Section \ref{s4}, we end with a few concluding remarks.\\

\section{Main Results}\label{s2}
In this section, we state our main results.\\
\begin{theorem}\label{thm21}
For every positive integer $L$, we have\\
\begin{align*}
\sum\limits_{n\ge 0}a_2(L,n)q^n = \dfrac{1}{(q;q^2)_L}\left(q^2+q^3\dfrac{1-q^{2L-2}}{1-q^2}+q^6\dfrac{1-q^{4L-4}}{1-q^4}\right).\numberthis\label{eq21}\\    
\end{align*}
\end{theorem}
\begin{remark}\label{rmk1}
As $L\rightarrow\infty$ in \eqref{eq21}, we get \eqref{eq13}.\\   
\end{remark}
\begin{theorem}\label{thm22}
For every positive integer $L$, we have\\
\begin{align*}
\sum\limits_{n\ge 0}b_2(L,n)q^n = q^2\dfrac{1-q^{L-1}}{1-q}(-q^2;q)_{L-2}.\numberthis\label{eq22}\\    
\end{align*}
\end{theorem}
\begin{remark}\label{rmk2}
As $L\rightarrow\infty$ in \eqref{eq22}, we get \eqref{eq14}.\\   
\end{remark}
\begin{theorem}\label{thm23}
For every positive integer $L$, $a_2(L,n)\ge b_2(L,n)$ for all $n\ge 0$.\\
\end{theorem}
Now, if $m$ is the total number of hooks of length $2$, we count $a_2(L,m,n)$ and $b_2(L,m,n)$ with a weight of $\binom{m}{2}$ and compute the generating functions of $\displaystyle{\sum\limits_{m\ge 0}\binom{m}{2}a_2(L,m,n)}$ and $\displaystyle{\sum\limits_{m\ge 0}\binom{m}{2}b_2(L,m,n)}$ respectively.\\
\begin{theorem}\label{thm24}
For every positive integer $L$, we have\\
\begin{align*}
\sum\limits_{n,m\ge 0}\binom{m}{2}a_2(L,m,n)q^n &= \dfrac{1}{(q;q^2)_L}\Bigg(q^5\dfrac{1-q^{2L-2}}{1-q^2}+(q^6+q^8)\dfrac{1-q^{4L-4}}{1-q^4}+q^8\left[\begin{matrix}L-1\\2\end{matrix}\right]_{q^2}\\&\quad\quad + q^{11}\dfrac{1-q^2}{1-q^{2L}}\left[\begin{matrix}L\\3\end{matrix}\right]_{q^2}((1+q^2)(1+q^{2L-2})+q^2+q^{2L})\\&\quad\quad + q^{16}\left[\begin{matrix}L-1\\2\end{matrix}\right]_{q^4}\Bigg).\numberthis\label{eq23}\\   
\end{align*}
\end{theorem}
\begin{remark}\label{rmk3}
It can be shown that the term on the right-hand side of \eqref{eq23}\\
\begin{align*}
\dfrac{1-q^2}{1-q^{2L}}\left[\begin{matrix}L\\3\end{matrix}\right]_{q^2}((1+q^2)(1+q^{2L-2})+q^2+q^{2L})\\    
\end{align*}
is non-negative. We leave it as an exercise for the motivated reader.\\
\end{remark}
\begin{theorem}\label{thm25}
For every positive integer $L$, we have\\
\begin{align*}
\sum\limits_{n,m\ge 0}\binom{m}{2}b_2(L,m,n)q^n = q^6\left[\begin{matrix}L-2\\2\end{matrix}\right]_{q}(-q^3;q)_{L-4}.\numberthis\label{eq24}\\ 
\end{align*}
\end{theorem}
\par We now state a conjecture showing bias of $\displaystyle{\sum\limits_{m\ge 0}\binom{m}{2}a_2(L,m,n)}$ over $\displaystyle{\sum\limits_{m\ge 0}\binom{m}{2}b_2(L,m,n)}$ for all $n\ge 0$.\\
\begin{conj}\label{conj26}
For every positive integer $L$, we have\\
\begin{align*}
\sum\limits_{m\ge 0}\binom{m}{2}a_2(L,m,n)\ge \sum\limits_{m\ge 0}\binom{m}{2}b_2(L,m,n),\\    
\end{align*}
for all $n\ge 0$.\\
\end{conj}
\begin{cor}\label{cor27}
We have\\
\begin{align*}
\sum\limits_{n,m\ge 0}\binom{m}{2}a_2(m,n)q^n &= \dfrac{1}{(q;q^2)_{\infty}}\Bigg(\dfrac{q^5}{1-q^2}+\dfrac{q^6+q^8}{1-q^4}+\dfrac{q^8}{(1-q^2)(1-q^4)}\\&\quad\quad + \dfrac{q^{11}+2q^{13}}{(1-q^4)(1-q^6)}+\dfrac{q^{16}}{(1-q^4)(1-q^8)}\Bigg).\numberthis\label{eq25}\\   
\end{align*}
\end{cor}
\begin{cor}\label{cor28}
We have\\
\begin{align*}
\sum\limits_{n,m\ge 0}\binom{m}{2}b_2(m,n)q^n = \dfrac{q^6}{(q;q^2)_{\infty}(1-q^2)(1-q^4)}.\numberthis\label{eq26}\\ 
\end{align*}
\end{cor}
\begin{theorem}\label{thm29}
For every non-negative integer $n$, we have\\
\begin{align*}
\sum\limits_{m\ge 0}\binom{m}{2}a_2(m,n)\ge \sum\limits_{m\ge 0}\binom{m}{2}b_2(m,n),\\
\end{align*}
for all $n\ge 0$.\\
\end{theorem}

\section{Proofs of Main Results}\label{s3}
In this section, we present our proofs of the main results.\\
\subsection{Proof of Theorem \ref{thm21}}\label{ss31}
Clearly,\\
\begin{align*}
a_2(L,n) = \sum\limits_{m\ge 0}ma_2(L,m,n).\\
\end{align*}
\par If $\pi_o$ is a partition into odd parts, then the number of hooks of length $2$ is equal to the number of different part sizes of $\pi_o$ that are greater than $1$ plus the number of different part sizes of $\pi_o$ that occur atleast twice. Let $F_2(L,z;q)$ denote the bi-variate generating function of $a_2(L,m,n)$. Then, we have\\
\begin{align*}
F_2(L,z;q) &:= \sum\limits_{n,m\ge 0}a_2(L,m,n)z^mq^n\\
&= \left(1+q+\dfrac{zq^2}{1-q}\right)\prod\limits_{n=2}^{L}\left(1+zq^{2n-1}+\dfrac{z^2q^{2(2n-1)}}{1-q^{2n-1}}\right),\numberthis\label{eq31}\\
\end{align*}
where \eqref{eq31} follows from standard combinatorial arguments in \cite{A98}.\\\par From the definition of $a_2(L,m,n)$, we have\\
\begin{align*}
\sum\limits_{n\ge 0}a_2(L,n)q^n = \dfrac{\partial}{\partial z}\at[\bigg]{z=1}F_2(L,z;q).\\
\end{align*}
Using logarithmic differentiation on both sides of \eqref{eq31} with respect to $z$ and setting $z=1$, we get\\
\begin{align*}
\sum\limits_{n\ge 0}a_2(L,n)q^n &= \dfrac{1}{(q;q^2)_L}\left(q^2+\sum\limits_{n=2}^{L}\left(q^{2n-1}+q^{4n-2}\right)\right)\\
&= \dfrac{1}{(q;q^2)_L}\left(q^2+q^3\dfrac{1-q^{2L-2}}{1-q^2}+q^6\dfrac{1-q^{4L-4}}{1-q^4}\right)\numberthis\label{eq32}\\
\end{align*}
where the last line of $\eqref{eq32}$ follows using geometric series. Thus, we have a proof of Theorem \ref{thm21}.\\\qed\\

\subsection{Proof of Theorem \ref{thm22}}\label{ss32}
Clearly,\\
\begin{align*}
b_2(L,n) = \sum\limits_{m\ge 0}mb_2(L,m,n).\\
\end{align*}
\par In the Young diagram of a partition $\pi_d = (\lambda_1,\lambda_2,\ldots)$ into distinct parts, the number of hooks of length $2$ in $\pi_d$ is equal to the number of parts $\lambda_i$ such that $\lambda_i - \lambda_{i+1}\ge 2$, where $\lambda_k = 0$ if $k > \nu(\pi_d)$. Thus, the number of hooks of length $2$ in $\pi_d$ where $l(\pi_d)\le L$ can be calculated as follows: start with the Young diagram of $\pi_d$ and remove the \textit{Sylvester triangle} by subtracting $1$ from the last part, $2$ from the second last part, etc., to obtain an ordinary partition $\mu$ having largest part $l(\mu)\le L-n$, where $n$ is the number of rows of the Young diagram of $\pi_d$ from which $1,2,\ldots,n$ are removed. Then, the \textit{Sylvester triangle} obtained is of size equal to\\
\begin{align*}
1+2+\ldots+n = \dfrac{n(n+1)}{2}.\\
\end{align*}
Now, we count the number of different part sizes in $\mu$ which can easily be seen to be equal to the number of hooks of length $2$ in $\pi_d$.\\\par Let $G_2(L,z;q)$ denote the bi-variate generating function of $b_2(L,m,n)$, i.e.,\\
\begin{align*}
G_2(L,z;q) &:= \sum\limits_{n,m\ge 0}b_2(L,m,n)z^mq^n.\\    
\end{align*}
\par Then, we have the following proposition\\
\begin{prop}\label{prop31}
\begin{align*}
G_2(L,z;q) = \sum\limits_{n=1}^{L}q^{\frac{n(n+1)}{2}}\sum\limits_{j=0}^{n}q^{\frac{j(j+1)}{2}}(z-1)^j\left[\begin{matrix}L-j\\n\end{matrix}\right]_{q}\left[\begin{matrix}n\\j\end{matrix}\right]_{q}.\numberthis\label{eq33}\\   
\end{align*}
\end{prop}
To prove Proposition \ref{prop31}, we need the following lemma.\\
\begin{lemma} (Alladi-Berkovich \cite[Lemma $3$]{AB04})\label{lemma31}
Let $h_{L,i}(y,q)$ denote the generating function of unrestricted partitions $\pi$ having $l(\pi)\le L$ into exactly $i$ non-negative parts such that $\pi$ is counted with weight $(1-y)^{\nu^{+}_d(\pi)}$, where $\nu^{+}_d(\pi)$ is the number of different positive parts of $\pi$. Then\\
\begin{align*}
h_{L,i}(y,q) = \sum\limits_{j\ge 0}q^{\frac{j(j+1)}{2}}(-y)^j\left[\begin{matrix}L+i-j\\i\end{matrix}\right]_{q}\left[\begin{matrix}i\\j\end{matrix}\right]_{q}.\numberthis\label{eq34}\\
\end{align*}
\end{lemma}
We now return to the proof of Proposition \ref{prop31}.\\\par
\textit{Proof}: It is easy to see that\\
\begin{align*}
G_2(L,z;q) = \sum\limits_{n=1}^{L}h_{L-n,n}(1-z,q)q^{\frac{n(n+1)}{2}},\\    
\end{align*}
which together with \eqref{eq34} yields \eqref{eq33}.\\\qed\\
\\\par From the definition of $b_2(L,m,n)$, we have\\
\begin{align*}
\sum\limits_{n\ge 0}b_2(L,n)q^n = \dfrac{\partial}{\partial z}\at[\bigg]{z=1}G_2(L,z;q).\\
\end{align*}    
Thus, from Proposition \ref{prop31}, we have\\
\begin{align*}
\sum\limits_{n\ge 0}b_2(L,n)q^n &= \sum\limits_{n=1}^{L}q^{\frac{n(n+1)}{2}+1}\left[\begin{matrix}L-1\\n\end{matrix}\right]_{q}\left[\begin{matrix}n\\1\end{matrix}\right]_{q}\\
&= q\dfrac{1-q^{L-1}}{1-q}\sum\limits_{n=1}^{L}q^{\frac{n(n+1)}{2}}\left[\begin{matrix}L-2\\n-1\end{matrix}\right]_{q}\\
&= q^2\dfrac{1-q^{L-1}}{1-q}\sum\limits_{n=0}^{L-2}q^{\frac{n(n+3)}{2}}\left[\begin{matrix}L-2\\n\end{matrix}\right]_{q}\\
&= q^2\dfrac{1-q^{L-1}}{1-q}(-q^2;q)_{L-2},\numberthis\label{eq35}\\
\end{align*}
where $\eqref{eq35}$ above follows by substituting $(z,N)\longmapsto (-q^2,L-2)$ in \cite[p.~$36$, eq. ($3.3.6$)]{A98}. Thus, we have a proof of Theorem \ref{thm22}.\\\qed\\

\subsection{Proof of Theorem \ref{thm23}}\label{ss33}
From Theorems \ref{thm21} and \ref{thm22}, we have\\
\begin{align*}
\sum\limits_{n\ge 0}\left(a_2(L,n)-b_2(L,n)\right)q^n &= \dfrac{1}{(q;q^2)_L}\left(q^2+q^3\dfrac{1-q^{2L-2}}{1-q^2}+q^6\dfrac{1-q^{4L-4}}{1-q^4}\right)\\
&\qquad - q^2\dfrac{1-q^{L-1}}{1-q}(-q^2;q)_{L-2}\\
& = q^2\left(\dfrac{1}{(q;q^2)_L}\left(1+q\dfrac{1-q^{2L-2}}{1-q^2}\right)-\dfrac{1-q^{L-1}}{1-q}(-q^2;q)_{L-2}\right)\\
&\qquad + \dfrac{q^6(1-q^{4L-4})}{(1-q^4)(q;q^2)_L}\\
& = q^2(A_L(q)-B_L(q))+C_L(q),\numberthis\label{eq36}\\
\end{align*}
where
\begin{align*}
A_L(q) = \dfrac{1}{(q;q^2)_L}\left(1+q\dfrac{1-q^{2L-2}}{1-q^2}\right),\quad B_L(q) = (-q^2;q)_{L-2}\dfrac{1-q^{L-1}}{1-q},     
\end{align*}
and
\begin{align*}
C_L(q) = \dfrac{q^6(1-q^{4L-4})}{(1-q^4)(q;q^2)_L}.\\    
\end{align*}
\par It is clear that $C_L(q)$, when expanded as a $q$-series, has non-negative coefficients, i.e., $C_L(q)\succeq 0$.
\\\par Now, define $\mathcal{A}_{L,n}$ to be the set of all pair partitions $(\pi_1,\pi_2)$ of $n$ such that $|\pi_1|+|\pi_2|=n$ where $\pi_1$ is a partition into blue odd parts having $l(\pi_1)\le 2L-1$ and $\pi_2$ is a partition into red odd parts having $l(\pi_2)\le 2L-3$ and $\nu(\pi_2)\le 1$.
\\\par Similarly, define $\mathcal{B}_{L,n}$ to be the set of all pair partitions $(\pi_1,\pi_2)$ of $n$ such that $|\pi_1|+|\pi_2|=n$ where $\pi_1$ is a partition into blue distinct parts different from $1$ where $l(\pi_1)\le L-1$ and $\pi_2$ is a partition into red parts having $l(\pi_2)\le L-2$ and $\nu(\pi_2)\le 1$. 
\\\par Note that $A_L(q)$ (resp. $B_L(q)$) is the generating function for pair partitions $(\pi_1,\pi_2)\in \mathcal{A}_{L,n}$ (resp. $\mathcal{B}_{L,n}$). We will now show that $A_L(q)\succeq B_L(q)$ by providing an explicit injection from $\mathcal{B}_{L,n}$ to $\mathcal{A}_{L,n}$.
\\\par As a corollary, we will have the desired inequality in Theorem \ref{thm23}.
\\\par For every non-negative integer $n$, we consider the following map:\\
\begin{align*}
\varphi = \varphi_{L,n} : \mathcal{B}_{L,n}\longrightarrow \mathcal{A}_{L,n}    
\end{align*}\\
Let $(\pi_1,\pi_2)\in\mathcal{B}_{L,n}$ and define $\varphi$ according to the following cases:
\\\par\begin{enumerate}
    \item \fbox{Case I: $\pi_2 = \emptyset$}\\\par In this case, define\\
    \begin{align*}
        \varphi((\pi_1,\pi_2)) = \varphi((\pi_1,\emptyset)) := (\sigma(\pi_1),\emptyset),\\
    \end{align*}
    where $\sigma$ is Sylvester's map as defined in section \ref{s1}.\\\par Clearly, $\varphi((\pi_1,\pi_2))\in\mathcal{A}_{L,n}$ because $l(\sigma(\pi_1))\le 2l(\pi_1)-1\le 2(L-1)-1 = 2L-3\le 2L-1$ where $l(\pi_1)\le L-1$, $l(\emptyset) = 0\le 2L-3$, $\nu(\emptyset) = 0\le 1$, and $|\sigma(\pi_1)|+|\emptyset| = |\pi_1|+0 = n$.\\
    \item \fbox{Case II: $\pi_2 = (\lambda_1)$ where $\lambda_1$ $(\ge 1)$ is odd}\\\par In this case, define\\
    \begin{align*}
        \varphi((\pi_1,\pi_2)) = \varphi((\pi_1,(\lambda_1))) := (\sigma(\pi_1),(\lambda_1)),\\
    \end{align*}
    where $\sigma$ is Sylvester's map as defined in section \ref{s1}.\\\par Clearly, $\varphi((\pi_1,\pi_2))\in\mathcal{A}_{L,n}$ because $l(\sigma(\pi_1))\le 2l(\pi_1)-1\le 2(L-1)-1 = 2L-3\le 2L-1$ where $l(\pi_1)\le L-1$, $\lambda_1$ is odd, $l((\lambda_1)) = \lambda_1\le L-2\le 2L-3$, $\nu((\lambda_1)) = 1\le 1$, and $|\sigma(\pi_1)|+|(\lambda_1)| = |\pi_1|+\lambda_1 = n$.\\
    \item \fbox{Case III: $\pi_2 = (\lambda_1)$ where $\lambda_1$ $(\ge 2)$ is even}\\\par In this case,\\
    \begin{align*}
        \varphi((\pi_1,\pi_2)) = \varphi((\pi_1,(\lambda_1))) := (\sigma(\pi^{*}_1),(\lambda_1-1)),\\
    \end{align*}
    where $\sigma$ is Sylvester's map as defined in section \ref{s1}, $\pi^{*}_1$ has all parts of $\pi_1$ which are different from $1$ and in addition, $\pi^{*}_1$ has the unique smallest part equal to $1$.\\\par Clearly, $\varphi((\pi_1,\pi_2))\in\mathcal{A}_{L,n}$ because $l(\sigma(\pi^{*}_1))\le 2l(\pi_1)-1\le 2(L-1)-1 = 2L-3\le 2L-1$ where $l(\pi^{*}_1)\le L-1$, $\lambda_1-1$ is odd, $l((\lambda_1-1)) = \lambda_1-1\le L-2\le 2L-3$, $\nu((\lambda_1-1)) = 1\le 1$, and $|\sigma(\pi^{*}_1)|+|(\lambda_1-1)| = |\pi_1|+1+\lambda_1-1 = n$.\\
\end{enumerate}
\par Then, it is easy to verify the following lemma\\
\begin{lemma}\label{lem32}
For every non-negative integer $n$,\\
\begin{align*}
\varphi = \varphi_{L,n} : \mathcal{B}_{L,n}\longrightarrow\mathcal{A}_{L,n}\\    
\end{align*}
is an injection.\\
\end{lemma}
\textit{Proof}:
It is easy to see that $\varphi$ is an injection. Let $(\pi_1,\pi_2)$, $(\Tilde{\pi}_{1},\Tilde{\pi}_{2})\in\mathcal{B}_{L,n}$. Then\\
\begin{align*}
\varphi((\pi_1,\pi_2)) = \varphi((\Tilde{\pi}_{1},\Tilde{\pi}_{2}))    
\end{align*}
implies
\begin{align*}
(\pi_1,\pi_2) = (\Tilde{\pi}_{1},\Tilde{\pi}_{2}).\\   
\end{align*}\qed\\
\par We now illustrate $\varphi$ with three examples:\\
\begin{example}
Let $L = 7$ and $n = 11$. Consider $(\pi_1,\pi_2) = ((\color{blue}6\color{black},\color{blue}5\color{black}), \color{red}\emptyset\color{black})\in\mathcal{B}_{7,11}$. Then, by Case I, we have\\
\begin{align*}
\varphi((\pi_1,\pi_2)) &= \varphi(((\color{blue}6\color{black},\color{blue}5\color{black}),\color{red}\emptyset\color{black}))\\
&= (\sigma((\color{blue}6\color{black},\color{blue}5\color{black})),\color{red}\emptyset\color{black})\\
&= ((\color{blue}{11}\color{black}),\color{red}\emptyset\color{black})\in\mathcal{A}_{7,11}.
\end{align*}
\begin{figure}[H]
$\ytableausetup{centertableaux}
\left(\,\,\ytableaushort
{{}{}{}{}{}{},{}{}{}{}{}}
* {6,5}
* [*(blue)]{6,5}
\,\,,\,\,
\color{red}\emptyset\color{black}\,\,\right)
\qquad\xlongrightarrow{\varphi}$
\\\quad\\\quad\\
$\left(\,\,\ytableaushort
{{}{}{}{}{}{}{}{}{}{}{}}
* {11}
* [*(blue)]{11}
\,\,,\,\,
\color{red}\emptyset\color{black}\,\,\right)$
\par\quad\\
\caption{$\varphi(((\color{blue}6\color{black},\color{blue}5\color{black}),\color{red}\emptyset\color{black})) = ((\color{blue}11\color{black}),\color{red}\emptyset\color{black})$}
\label{fig5}
\end{figure}
\end{example}\quad\\
\begin{example}
Let $L = 7$ and $n = 16$. Consider $(\pi_1,\pi_2) = ((\color{blue}6\color{black},\color{blue}5\color{black}), (\color{red}5\color{black}))\in\mathcal{B}_{7,16}$. Then, by Case II, we have\\
\begin{align*}
\varphi((\pi_1,\pi_2)) &= \varphi(((\color{blue}6\color{black},\color{blue}5\color{black}), (\color{red}5\color{black})))\\
&= (\sigma((\color{blue}6\color{black},\color{blue}5\color{black})), (\color{red}5\color{black}))\\
&= ((\color{blue}{11}\color{black}), (\color{red}5\color{black}))\in\mathcal{A}_{7,16}.
\end{align*}
\begin{figure}[H]
$\ytableausetup{centertableaux}
\left(\,\,\ytableaushort
{{}{}{}{}{}{},{}{}{}{}{}}
* {6,5}
* [*(blue)]{6,5}
\,\,,\,\,
\ytableaushort
{{}{}{}{}{}}
* {5}
* [*(red)]{5}\,\,\right)
\qquad\xlongrightarrow{\varphi}$
\\\quad\\\quad\\
$\left(\,\,\ytableaushort
{{}{}{}{}{}{}{}{}{}{}{}}
* {11}
* [*(blue)]{11}
\,\,,\,\,
\ytableaushort
{{}{}{}{}{}}
* {5}
* [*(red)]{5}\,\,\right)$
\par\quad\\
\caption{$\varphi(((\color{blue}6\color{black},\color{blue}5\color{black}),(\color{red}5\color{black}))) = ((\color{blue}11\color{black}),(\color{red}5\color{black}))$}
\label{fig6}
\end{figure}
\end{example}\quad\\
\begin{example}
Let $L = 7$ and $n = 13$. Consider $(\pi_1,\pi_2) = ((\color{blue}6\color{black},\color{blue}5\color{black}), (\color{red}2\color{black}))\in\mathcal{B}_{7,13}$. Then, by Case III, we have\\
\begin{align*}
\varphi((\pi_1,\pi_2)) &= \varphi(((\color{blue}6\color{black},\color{blue}5\color{black}), (\color{red}2\color{black})))\\
&= (\sigma((\color{blue}6\color{black},\color{blue}5\color{black},\color{blue}1\color{black})), (\color{red}1\color{black}))\\
&= ((\color{blue}9\color{black},\color{blue}3\color{black}), (\color{red}1\color{black}))\in\mathcal{A}_{7,13}.
\end{align*}
\begin{figure}[H]
$\ytableausetup{centertableaux}
\left(\,\,\ytableaushort
{{}{}{}{}{}{},{}{}{}{}{}}
* {6,5}
* [*(blue)]{6,5}
\,\,,\,\,
\ytableaushort
{{}{}}
* {2}
* [*(red)]{2}\,\,\right)
\qquad\xlongrightarrow{\varphi}$
\\\quad\\\quad\\
$\left(\,\,\ytableaushort
{{}{}{}{}{}{}{}{}{},{}{}{}}
* {9,3}
* [*(blue)]{9,3}
\,\,,\,\,
\ytableaushort
{{}}
* {1}
* [*(red)]{1}\,\,\right)$
\par\quad\\
\caption{$\varphi(((\color{blue}6\color{black},\color{blue}5\color{black}),(\color{red}2\color{black}))) = ((\color{blue}9\color{black},\color{blue}3\color{black}),(\color{red}1\color{black}))$}
\label{fig7}
\end{figure}
\end{example}\quad\\
Thus, Lemma \ref{lem32} implies the inequality in Theorem \ref{thm23}.\\\qed\\

\subsection{Proof of Theorem \ref{thm24}}\label{ss34}
From the definition of $F_2(L,z;q)$ in \eqref{eq31}, observe that\\
\begin{align*}
\sum\limits_{n\ge 0}\left(\sum\limits_{m\ge 0}(m^2-m)a_2(L,m,n)\right)q^n = \dfrac{\partial^2}{\partial z^2}\at[\bigg]{z=1}F_2(L,z;q).\\    
\end{align*}
Using logarithmic differentiation on both sides of \eqref{eq31} with respect to $z$, we have\\
\begin{align*}
\dfrac{\partial}{\partial z}F_2(L,z;q) = F_2(L,z;q)\left(\dfrac{q^2}{1+(z-1)q^2}+\sum\limits_{n=2}^{L}\dfrac{q^{2n-1}+(2z-1)q^{4n-2}}{1+(z-1)q^{2n-1}+(z^2-z)q^{4n-2}}\right).\numberthis\label{eq37}\\    
\end{align*}
Now, applying logarithmic differentiation on both sides of \eqref{eq37} with respect to $z$ and setting $z=1$, we have\\
\begin{align*}
\sum\limits_{n\ge 0}\left(\sum\limits_{m\ge 0}(m^2-m)a_2(L,m,n)\right)q^n &= \dfrac{1}{(q;q^2)_L}\left(q^2+q^3\dfrac{1-q^{2L-2}}{1-q^2}+q^6\dfrac{1-q^{4L-4}}{1-q^4}\right)^2\numberthis\label{eq38}\\
&\quad + \dfrac{1}{(q;q^2)_L}\Bigg(-q^4+q^6\dfrac{1-q^{4L-4}}{1-q^4}-2q^9\dfrac{1-q^{6L-6}}{1-q^6}\\
&\quad - q^{12}\dfrac{1-q^{8L-8}}{1-q^8}\Bigg).\\
\end{align*}
Dividing both sides of \eqref{eq38} by $2$ and simplifying, we get\\
\begin{align*}
\sum\limits_{n\ge 0}\left(\sum\limits_{m\ge 0}\binom{m}{2}a_2(L,m,n)\right)q^n &= \dfrac{1}{(q;q^2)_L}\Bigg(q^5\dfrac{1-q^{2L-2}}{1-q^2}+(q^6+q^8)\dfrac{1-q^{4L-4}}{1-q^4}\\
&\quad + q^8\left[\begin{matrix}L-1\\2\end{matrix}\right]_{q^2}+q^{11}\dfrac{1-q^2}{1-q^{2L}}\left[\begin{matrix}L\\3\end{matrix}\right]_{q^2}((1+q^2)(1+q^{2L-2})\\
&\quad + q^2+q^{2L})+q^{16}\left[\begin{matrix}L-1\\2\end{matrix}\right]_{q^4}\Bigg).\\    
\end{align*}
Thus, we have a proof of Theorem \ref{thm24}.\\\qed\\

\subsection{Proof of Theorem \ref{thm25}}\label{ss35}
From the definition of $G_2(L,z;q)$ in \eqref{eq33}, observe that\\
\begin{align*}
\sum\limits_{n\ge 0}\left(\sum\limits_{m\ge 0}(m^2-m)b_2(L,m,n)\right)q^n = \dfrac{\partial^2}{\partial z^2}\at[\bigg]{z=1}G_2(L,z;q).\\    
\end{align*}
Differentiating both sides of \eqref{eq33} with respect to $z$, we have\\
\begin{align*}
\dfrac{\partial}{\partial z}G_2(L,z;q) = \sum\limits_{n=1}^{L}q^{\frac{n(n+1)}{2}}\sum\limits_{j=1}^{n}jq^{\frac{j(j+1)}{2}}(z-1)^{j-1}\left[\begin{matrix}L-j\\n\end{matrix}\right]_{q}\left[\begin{matrix}n\\j\end{matrix}\right]_{q}.\numberthis\label{eq39}\\    
\end{align*}
Now, differentiating both sides of \eqref{eq39} with respect to $z$ and setting $z=1$, we have\\
\begin{align*}
\sum\limits_{n\ge 0}\left(\sum\limits_{m\ge 0}(m^2-m)b_2(L,m,n)\right)q^n = 2\sum\limits_{n=2}^{L}q^{\frac{n(n+1)}{2}+3}\left[\begin{matrix}L-2\\n\end{matrix}\right]_{q}\left[\begin{matrix}n\\2\end{matrix}\right]_{q}.\numberthis\label{eq10}\\    
\end{align*}
Dividing both sides of \eqref{eq10} by $2$, we get\\
\begin{align*}
\sum\limits_{n\ge 0}\left(\sum\limits_{m\ge 0}\binom{m}{2}b_2(L,m,n)\right)q^n &= \sum\limits_{n=2}^{L}q^{\frac{n(n+1)}{2}+3}\left[\begin{matrix}L-2\\n\end{matrix}\right]_{q}\left[\begin{matrix}n\\2\end{matrix}\right]_{q}\\
&= q^3\dfrac{(1-q^{L-3})(1-q^{L-2})}{(1-q)(1-q^2)}\sum\limits_{n=2}^{L}q^{\frac{n(n+1)}{2}}\left[\begin{matrix}L-4\\n-2\end{matrix}\right]_{q}\\
&= q^6\left[\begin{matrix}L-2\\2\end{matrix}\right]_{q}\sum\limits_{n=0}^{L-2}q^{\frac{n^2+5n}{2}}\left[\begin{matrix}L-4\\n\end{matrix}\right]_{q}\\
&= q^6\left[\begin{matrix}L-2\\2\end{matrix}\right]_{q}\sum\limits_{n=0}^{L-4}q^{\frac{n^2+5n}{2}}\left[\begin{matrix}L-4\\n\end{matrix}\right]_{q}\\
&= q^6\left[\begin{matrix}L-2\\2\end{matrix}\right]_{q}(-q^3;q)_{L-4},\numberthis\label{eq311}\\
\end{align*}
where \eqref{eq311} above follows by substituting $(z,N)\longmapsto (-q^3,L-4)$ in \cite[p.~$36$, eq. ($3.3.6$)]{A98}. Thus, we have a proof of Theorem \ref{thm25}.\\\qed\\

\subsection{Proof of Corollary \ref{cor27}}\label{ss36}
Letting $L\rightarrow\infty$ on both sides of \eqref{eq23}, we get\\
\begin{align*}
\sum\limits_{n,m\ge 0}\binom{m}{2}a_2(m,n)q^n &= \dfrac{1}{(q;q^2)_{\infty}}\Bigg(\dfrac{q^5}{1-q^2}+\dfrac{q^6+q^8}{1-q^4}+\dfrac{q^8}{(1-q^2)(1-q^4)}\\&\quad\quad + \dfrac{q^{11}+2q^{13}}{(1-q^4)(1-q^6)}+\dfrac{q^{16}}{(1-q^4)(1-q^8)}\Bigg),\numberthis\label{eq312}\\
\end{align*}
where we use the fact that\\
\begin{align*}
\lim_{n\rightarrow\infty}\left[\begin{matrix}n\\m\end{matrix}\right]_{q} = \dfrac{1}{(q;q)_m}.\\    
\end{align*}
Thus, we have a proof of Corollary \ref{cor27}.\\\qed\\

\subsection{Proof of Corollary \ref{cor28}}\label{ss37}
Letting $L\rightarrow\infty$ on both sides of \eqref{eq24}, we get\\
\begin{align*}
\sum\limits_{n,m\ge 0}\binom{m}{2}a_2(m,n)q^n &= \dfrac{q^6(-q^3;q)_{\infty}}{(1-q)(1-q^2)}\\
&= \dfrac{q^6(-q;q)_{\infty}}{(1-q^2)(1-q^4)}\\
&= \dfrac{q^6}{(q;q^2)_{\infty}(1-q^2)(1-q^4)},\numberthis\label{eq313}\\
\end{align*}
where we use the fact that\\
\begin{align*}
\lim_{n\rightarrow\infty}\left[\begin{matrix}n\\m\end{matrix}\right]_{q} = \dfrac{1}{(q;q)_m}\\    
\end{align*}
and \eqref{eq313} follows using Euler's identity \cite[p.~$5$, Corollary $1.2$]{A98}. Thus, we have a proof of Corollary \ref{cor28}.\\\qed\\

\subsection{Proof of Theorem \ref{thm29}}\label{ss38}
From Corollaries \ref{cor27} and \ref{cor28}, we have\\
\begin{align*}
\sum\limits_{n,m\ge 0}\binom{m}{2}&(a_2(m,n)-b_2(m,n))q^n = \dfrac{1}{(q;q^2)_{\infty}}\Bigg(\dfrac{q^5}{1-q^2}+\dfrac{q^6+q^8}{1-q^4}+\dfrac{q^8}{(1-q^2)(1-q^4)}\\&\qquad\qquad\qquad\qquad\qquad + \dfrac{q^{11}+2q^{13}}{(1-q^4)(1-q^6)}+\dfrac{q^{16}}{(1-q^4)(1-q^8)}\Bigg)\\
&\qquad\qquad\qquad\qquad\qquad - \dfrac{q^6}{(q;q^2)_{\infty}(1-q^2)(1-q^4)}\\
&= \dfrac{q^5(1-q^2)^2(1+2q^2+q^3+2q^4+q^5+2q^6+q^7+2q^8+q^{10}+q^{12}+q^{14})}{(q;q^2)_{\infty}(1-q^2)(1-q^4)(1-q^6)(1-q^8)}\\
&= \dfrac{q^5(1+q)(1+2q^2+q^3+2q^4+q^5+2q^6+q^7+2q^8+q^{10}+q^{12}+q^{14})}{(q^3;q^2)_{\infty}(1-q^4)(1-q^6)(1-q^8)}.\numberthis\label{eq314}\\
\end{align*}
Clearly, the expression in \eqref{eq314}, expanded as a $q$-series, has non-negative coefficients. Thus, we have a proof of Theorem \ref{thm29}.\\\qed\\

\section{Concluding Remarks}\label{s4}
It should be noted that Conjecture \ref{conj26} follows from the simpler inequality which can be stated as follows.\\
\begin{conj}\label{conj41}
For every positive integer $L$, we have\\
\begin{align*}
\dfrac{1}{(q;q^2)_L}\left(\dfrac{1-q^{2L-2}}{1-q^2}+q^3\left[\begin{matrix}L-2\\2\end{matrix}\right]_{q^2}\right)\succeq q\left[\begin{matrix}L-2\\2\end{matrix}\right]_{q}(-q^3;q)_{L-4}.\numberthis\label{eq41}\\    
\end{align*}
\end{conj}
Using Remark \ref{rmk3}, we see that Conjecture \ref{conj41} implies Conjecture \ref{conj26}.\\\par It would be interesting to find an injective proof of Conjecture \ref{conj41}.\\\par We plan to investigate finite analogs of \eqref{eq15} and \eqref{eq16} in a forthcoming paper and also check whether there is a bias between $a_3(L,n)$ and $b_3(L,n)$.\\

\section*{Acknowledgments}
We would like to thank Krishnaswami Alladi, George E. Andrews, Ali K. Uncu and Hamza Yesilyurt for their kind interest. We also thank Will Craig for directing us to \cite{CDH23} for a proof of Conjecture \ref{conj14}.\\

\bibliographystyle{amsplain}


\end{document}